\documentclass{article}

\usepackage{amsmath}
\usepackage{amssymb}
\usepackage{amsthm}
\usepackage{graphics}
\usepackage[latin1]{inputenc}
\usepackage[english]{babel}
\usepackage[T1]{fontenc}

\usepackage{xcolor}
\usepackage{comment}

\newtheorem{thm}{Theorem}[section]
\newtheorem{lemma}[thm]{Lemma}
\newtheorem{rem}[thm]{Remark}

\makeatletter

\newcommand*{\plim}[1][]{%
	\if\relax\detokenize{#1}\relax
	\def\next{\qopname\relax m{lim}}%
	\else
	\def\next{\qopname\newmcodes@ m{#1-lim}}%
	\fi
	\next
}

\newcommand*{\psum}[1][]{%
	\DOTSB
	\if\relax\detokenize{#1}\relax\else
	\operatorname{#1-}\mkern-\thinmuskip
	\fi
	\sum@\slimits@
}

\makeatother

\newcommand{\R}{\mathbb{R}}             
\newcommand{\N}{\mathbb{N}}             
\newcommand{\Z}{\mathbb{Z}}             
\newcommand{\Q}{\mathbb{Q}}             
\newcommand{\C}{\mathbb{C}}             

\newcommand{\Ab}{{\bf A}}

\newcommand{\half}{\frac{1}{2}}

\newcommand{\tr}{{\rm{Tr}}\ }

\newcommand {\clb}{\color{blue}}
\newcommand {\clr} {\color{red}}

\newcommand{\polylog}{\rm{Li}}
\newcommand{\m}{|m-\nu|}

\numberwithin{equation}{section}

\begin{document}
	
	\title{Trace formulas for the magnetic Laplacian and Dirichlet to Neumann operator -- Explicit expansions --}
	\author{Bernard Helffer  and 
	Fran\c cois Nicoleau \ \footnote{Research supported by the French National Research Project GDR Dynqua} \\[12pt]
    \small  Laboratoire de Math\'ematiques Jean Leray, UMR CNRS 6629, \\ 
		\small Nantes Universit\'e,  F-44000 Nantes, France.\\
		\small Email: bernard.helffer@univ-nantes.fr , francois.nicoleau@univ-nantes.fr \\
		}

	
	\date{\today}

\maketitle

\begin{abstract}
Inspired by a recent paper of G. Liu and X. Tan (2023), we would like to measure how the magnetic effect appears in the heat trace formula associated with the magnetic Laplacian and  the magnetic Dirichlet-to-Neumann operator. We propose to  the reader an overview of magnetic heat trace formulas through explicit examples. On the way we obtain new formulas and in particular we calculate explicitely some non local terms and logarithmic terms appearing in the Steklov heat trace asymptotics.


\vspace{0.5cm}




\end{abstract}


\tableofcontents


\section{Introduction}

In this paper, we study the heat trace asymptotics associated with the magnetic Steklov problem  on a  smooth compact Riemannian manifold $(\Omega,g)$ with $C^\infty$ boundary $\partial \Omega$,  (in particular this could be a bounded domain in $\mathbb R^d$).  We denote by \rm{grad}, \rm{div} and $\langle  \cdot , \cdot  \rangle$ the gradient operator, the divergence operator and the inner product on $\Omega$ with respect to the metric $g$ respectively. 

\vspace{0.2cm}
\noindent
For any $u \in C_0^{\infty}(\Omega)$, the magnetic Schr\"odinger operator is defined as 
\begin{equation}\label{defMagOp}
H_{A,V}\  u = -\Delta_g u -2i \  \langle A, \rm{grad} \ u \rangle + (A^2 -i \  \rm{div} \ A +V) u,
\end{equation}
where  $\Delta_g$ is the Laplace-Beltrami operator, ${\displaystyle{A= \sum_{j=1}^n A_j dx_j}}$ is the 1-form magnetic potential and $V$ is the electric potential. We often identify the $1$-form magnetic potential $A$ with the vector field 
$\overrightarrow{A} = (A_1, ..., A_n)$. The magnetic field is given by the $2$-form $B =dA$.

\vspace{0.2cm}
In the following, we assume that the real-valued potentials $A$ and $V$ are smooth in $\overline{\Omega}$. Nevertheless  we will discuss in some sections weaker assumptions\footnote{See   Leinfelder-Simader \cite{Le1983,LeSi1981}, where it is assumed 
$A \in L^2_{loc}$, $V$ in $L^2_{loc}$ and semi-bounded. See  also Simon, Subsection B.213 in \cite{Si1982}. }
and will focus in this direction on Aharonov-Bohm like singularities avoiding the boundary.

\vspace{0.2cm}
We also assume that zero does not belong to the spectrum of the Dirichlet realization  of  $H_{A,V}$, so that the boundary value  problem 
\begin{equation}  \label{Dirichlet}
	\left\{
	\begin{array}{rll}
		H_{A,V} \  u &=&0  \  \ \rm{in}  \ \ \Omega,\\
		u_{  \vert \partial \Omega} & =& f \in H^{1/2}(\partial\Omega) .
	\end{array}\right.
\end{equation}
has a unique solution $u \in H^1(\Omega)$. The  Dirichlet to Neumann map, (in what follows D-to-N map), is defined by
\begin{equation} \label{D-t-N--map}
\begin{array}{rll}
\Lambda_{A,V} :   H^{1/2}(\partial\Omega)& \longmapsto & H^{-1/2}(\partial\Omega) \\ 
	               f   &\longmapsto&  \left(\partial_{\nu} u + i \langle A, \nu \rangle \ u  \right)_{  \vert \partial \Omega} ,
\end{array}
\end{equation}
where $\nu$ is the outward normal unit vector field on $\partial\Omega$. 

\vspace{0.2cm}\noindent
Under a gauge transformation in the magnetic potential $A \to A + \rm{grad} \  \varphi$ with $e^{i \varphi} \in C^\infty (\overline{\Omega})$, 
one has
\begin{equation}
	\Lambda_{A + \rm{grad} \, \varphi  ,V}  =  e^{i \varphi} \Lambda_{A,V} e^{-i \varphi}.
\end{equation}
It follows that $\Lambda_{A,V}$ carries information about the magnetic field instead of information about the magnetic potential $A$. Morevover, 
(see for example Appendix in \cite{FoHe2010}), one can assume that   $\langle A, \nu \rangle =0$   on $\partial \Omega$. With this assumption the magnetic normal 
derivative $(\partial_{\nu}  +i  \langle A, \nu \rangle )$ becomes the standard normal derivative.

\vspace{0.2cm}
 
In the following, we are interested in recovering information about the magnetic field $B$ and the electric potential from the D-to-N map $\Lambda_{A,V}$ . This is a classical inverse problem where one wishes to know the interior properties of a medium by making  only electrical measurements at its boundary. It is related to the famous Calder\'on inverse conductivity problem \cite{Cal1980}. In \cite{NSU1995}, for smooth potentials $A$ and $V$,  G. Nakamura, Z. Sun and G. Uhlmann recover the electrical potential et the magnetic field from the D-to-N map. This smoothness assumption was reduced to $C^1$ by C. Tolmasky \cite{To1998} and to Dini continuous by M. Salo \cite{Sa2004}. We also mention the papers  \cite{BrSa2006,DSKSJ2007, Tz2008} which consider the boundary determination, partial Cauchy data and the stability for this magnetic inverse problem.

\vspace{0.2cm}
\noindent
We recall that the spectrum of the D-to-N is discrete and is given by a sequence of eigenvalues 
\begin{equation} \label{spectrum}
	\lambda_1 \leq \lambda_2 \leq ... \leq \lambda_n \leq ...  \to + \infty.
\end{equation}

\vspace{0.2cm}
As it was proved in \cite{LiuTan2023, NSU1995}, (see also   \cite{GKLP2022, PoSh2015}), when the Riemannian metric $g$, the electromagnetic potentials $A$, $V$ are smooth i.e. in $C^\infty(\overline{\Omega})$, the D-to-N map  $\Lambda_{A,V}$ is a self-adjoint elliptic classical pseudodifferential operator of order $1$ on the boundary, with the same principal symbol as $\sqrt{-\Delta_g} _{ | \partial \Omega}$, the square root of the  boundary Laplacian. We will discuss later how to relax the assumption of regularity. It follows from \cite{Hor1968} that the spectral counting function $N(\lambda)$ of the D-to-N map  satisfies the Weyl asymptotics :
 
\begin{equation}\label{counting}
N(\lambda) =  \frac{\rm{Vol} (B^{d-1} ) \ \rm{Vol}(\partial \Omega)}{(2\pi)^{d-1}}   \ \lambda^{d-1}  + \mathcal O (\lambda^{d-2})\ ,\ \lambda \to + \infty,
\end{equation}

which implies that the Steklov  eigenvalues $\lambda_n$ verify 
\begin{equation} \label{Weyl}
	\lambda_n = 2 \pi \left( \frac{n}{\rm{Vol} (B^{d-1} ) \ \rm{Vol}(\partial \Omega)}   \right)^{\frac{1}{d-1}} + \mathcal O(1) \ , 
	\ n \to + \infty. 
\end{equation}
Here, $B^{d-1}$ denotes the unit ball in $\R^{d-1}$. \\

Note that in the $(2D)$-case it is known that the D-to-N operator equals the square-root of the Laplace-Beltrami operator on the boundary modulo a regularizing operator. See  \cite{GKLP2022}.

\vspace{0.2cm}

As far as we know the dependence on the magnetic field in the Weyl's formula relative to the D-to-N operator has not been analyzed, (see nevertheless \cite{HaIv2017} for  accurate results in the D-to-N problem with non zero frequency).

\vspace{0.2cm}
In order to get information about the magnetic field and the electric potential, we consider the trace of the associated heat operator $e^{-t \Lambda_{A,V}}$  which 
admits the following asymptotic expansion as $t \to 0^+$, (see \cite{GiPo2017}, and references therein) : 
\begin{equation}
	\tr (e^{-t \Lambda_{A,V}} ) = \sum_{n=1}^{+\infty} e^{-t \lambda_n} \sim \sum_{k=0}^{\infty} a_k \ t^{-d+k+1} + \sum_{\ell =1}^{\infty} b_\ell \  t^\ell \log t.
\end{equation}
The coefficients $a_k$ and $b_\ell $ are called the Steklov heat invariants and can be theoretically computed using the so-called zeta function, (for simplicity, we assume here that $\Lambda_{A,V}$ is positive):  
\begin{equation}\label{zeta}
	\zeta( \Lambda_{A,V}, z)  := {\rm{Tr}} \left( \Lambda_{A,V}^{-z} \right) \ , \ {\rm{Re}} \ z \gg1.
\end{equation}
Thus, using the Mellin transform, one gets:
\begin{equation}
	\Gamma(z) \  \zeta( \Lambda_{A,V}, z)  = \int_0^{+\infty} t^{z-1} \ 	\tr  (e^{-t \Lambda_{A,V}} ) \ dt,
\end{equation}
and we can easily see that the function $$F(z) := \Gamma(z) \  \zeta( \Lambda_{A,V}, z)$$ has a meromorphic extension to $\C$ with possible  poles at the points $d-1-k$, with $ k \geq 0$ and at the points $-\ell $ with $\ell  \geq 1$.\\ Then, we have a relationship between the Steklov invariants and the residues of $F(z)$ at these possible poles:
\begin{equation}
	a_k = {\rm{Res}} \ ( F(z)\ ;\  d-1-k) \ ,\ b_\ell =  -{\rm{Res}} \ ( (z+\ell) F(z) \ ; \ -\ell )  .
\end{equation} 
Actually, the coefficients $a_k$ for $k=0, ... ,d-1$ and $b_\ell $  for $\ell  \geq 1$ are local invariants,  \cite{GiPo2017}. They are integral over the boundary  of  functions $a_k (x)$ and $b_\ell(x)$ which are  explicit polynomials   in the metric and its inverse, in the electromagnetic potentials and their derivatives in tangential and normal directions along the boundary, at the point $x$ :
\begin{equation}
	a_k=\int_{\partial \Omega} a_k(x) \ d\sigma \ ,\ b_\ell =\int_{\partial \Omega} b_\ell(x) \ d\sigma ,
\end{equation}
where $d \sigma$ is the induced measure on the boundary $\partial \Omega$.\\
\begin{rem}
As a consequence of the local character of  these $a_k$'s, the $a_k$'s are gauge invariants for $k=0,\dots,d-1$\,.
We conjecture that  there are  explicit polynomials   in the metric and its inverse, in the electromagnetic fields and their derivatives in tangential and normal directions along the boundary, at the point $x$. This should be a consequence of a magnetic pseudo-differential calculus extending Boutet de Monvel calculus (see \cite{Bou1971,MP2004,HP2010}). 
\end{rem}

In contrast, the coefficients $a_k$ for $k \geq d$ are not local invariants. We refer to \cite{GiGr1998, PoSh2015} for details. \\ 


The usual way to compute  the local heat  invariants $a_k$, $k=0, ... ,d-1$,  is to use the standard pseudo-differential calculus on manifolds applying the work\footnote{See also G. Grubb,  Israel J. Math. 1971  explaining
 the relation between the Calderon projector and the D-to-N operator.} of Seeley \cite{GiGr1998, Se1967}. The pseudo-differential operators on the boundary are classical (the symbol is asymptotically a sum of homogeneous symbols) and the symbolic calculus on the boundary  is always defined modulo regularizing operators, i.e. operators whose distribution kernel is $C^\infty$ on $\partial \Omega \times \partial \Omega$.  Since we are also interested with the non local invariants, we propose an alternative and  direct calculation.

To make it clearer, let us make the calculations with a toy model. As we shall see in Section~\ref{champconstant}, the Steklov eigenvalues $\lambda_n$ associated with a constant magnetic field of strenght $b$ in the unit disk of $\R^2$ are explicitely known. There are given in terms of the generalized Laguerre functions. In this case, the magnetic potential is given by $A = b\ (-ydx+xdy)$. The eigenvalues $\lambda_n$ satisfy the following asymptotics as $n \to + \infty$:
\begin{equation}
	\lambda_n = n +\alpha + \frac{\beta}{n} + \mathcal O(\frac{1}{n^2}),
\end{equation}
where $\alpha$ and $\beta$ are some suitable real constants. To simplify the notation, we set:
$$\mu_n=  n +\alpha + \frac{\beta}{n}\,.
$$ Then, the heat trace of the D-to-N map $\Lambda_{A,0}$ is given by 
\begin{eqnarray}\label{decompositiontrace}
	\tr (e^{-t \Lambda_{A,0}} ) &=& \sum_{n=1}^{+\infty} e^{-t \lambda_n} \nonumber \\
	  &=&  \sum_{n=1}^{+\infty} e^{-t \mu_n} + \sum_{n=1}^{+\infty} e^{-t \mu_n} \left( e^{-t (\lambda_n - \mu_n)}-1 \right).
\end{eqnarray}
The first term of the (RHS) of (\ref{decompositiontrace}) can be explicitely computed modulo $\mathcal O(t^3)$, (for instance).  One has:
\begin{eqnarray*}\label{toy}
\sum_{n=1}^{+\infty} e^{-t \mu_n}	&=& \sum_{n=1}^{+\infty} e^{-t(n+\alpha + \frac{\beta}{n})} = e^{-\alpha t} \sum_{n=1}^{+\infty} e^{-tn} \                                  e^{-\frac{t \beta }{n}}\\
	&=& e^{-\alpha t}  \sum_{n=1}^{+\infty} e^{-tn} \ \left( 1- \frac{\beta t}{n} + \frac{(\beta t)^2}{2n^2}+ \mathcal O(\frac{t^3}{n^3} ) \right) \\
	&=& e^{-\alpha t}  \left(  \frac{1}{e^{t}-1} - \beta t \ \polylog_1(e^{-t})  + \frac{(\beta t)^2}{2}
	\polylog_2(e^{-t}) \right) +  \mathcal O(t^3),
\end{eqnarray*}                                         
where $${\displaystyle{\polylog_s(z)= \sum_{n=1}^{+\infty} \frac{z^n}{n^s}}}$$ is the so-called polylogarithm  function.\\
 In particular,  one has $\polylog_1(z) =  - \log (1-z) $, and we get the following asymptotic expansion:
\begin{eqnarray}\label{asympttoy}
	\sum_{n=1}^{+\infty} e^{-t \mu_n}& =&  \frac{1}{t}  -\left( \alpha + \frac{1}{2}\right)  + \beta t \log t  + \frac{1}{12} \left(1+6 \alpha +6 \alpha^2 \right) t \nonumber \\
	&&	- \alpha \beta  \ t^2 \log t  -  \frac{1}{12} \left(\alpha+3 \alpha^2+2
	\alpha^3+6 \beta -\beta^2 \pi ^2\right) t^2  + \ \mathcal O(t^3).
\end{eqnarray} \\
The second term of the (RHS) of (\ref{decompositiontrace}) can be written as 
\begin{equation}
\sum_{n=1}^{+\infty} \left( \sum_{m=0}^{M }  \frac{(-t \mu_n)^m}{m!}  + O\left((t\mu_n)^{M+1}\right) \right) \cdot \left(  \sum_{p=1}^{P}  \frac{(-t (\lambda_n -\mu_n))^p}{p!} +  O \left( t (\lambda_n - \mu_n)^{P+1} \right) \right).
\end{equation}
If $\lambda_n - \mu_n$ decays sufficiently, we can invert  the order of summation giving us an additional  asymptotic expansion.  Actually, as we shall see in Section \ref{champconstant}, we can justify such an inversion up to $\mathcal O(t^3)$ if the asymptotics  of the Steklov eigenvalues  $\lambda_n$ are known modulo $\mathcal O(n^{-4})$. Of course, this affects the definition of $\mu_n$ and the asymptotics given in (\ref{asympttoy}). This leads to cumbersome calculations requiring the assistance of a computer.

\vspace{0.2cm}
The paper is organized as follows. In Section 2, we review for comparison the heat trace asymptotic expansions for several magnetic Hamiltonians, (including the singular case of the Aharonov-Bohm oscillator), first on the unit disk in $\R^2$, and secondly on the sphere $S^2$ in $\R^3$. In particular, Formulas (\ref{traceJacobi}) and (\ref{tracesphere}) seem to be new. In Section 3, we recall the general results and extend it to a more singular case.  In Sections 4-5, we calculate the  asymptotics of the heat trace operator for the D-to-N maps associated with Aharonov-Bohm Laplacians. In Section 6, we consider the case of a constant magnetic field in the unit disk in $\R^2$ and we calculate explicitely the first logarithmic terms appearing in the heat trace asymptotics. According to us, this last result result is also new.

\vspace{0.3cm}
\noindent \textbf{Acknowledgements}: We would like to warmly thank Gerd Grubb and El Maati Ouhabaz for enlightning exhanges of e-mails.\\

\section{On heat expansions for the magnetic Laplacian.}

\subsection{Diamagnetic inequality.}

In this section, we consider the magnetic Schrödinger operator $H_{A,V}$ defined in the introduction with $\Omega= \R^2$, (except for the subsection \ref{sectionsphere} where we consider the case of a closed Riemannian surface). We assume that the real-valued magnetic and electric potentials $A$ and $V$ are smooth. The diamagnetic inequality can be expressed as:
\begin{equation}\label{ineqmag}
	|\exp (-t H_{A,V}) f |\ (x) \leq( \exp (- t H_{0,V}) |f| ) \ (x)\,,
\end{equation}
for $t\geq 0$ and almost every $x$, (see for instance \cite{AHS1978, CSS1978, Ou1996, SiJFA1979}).

\vspace{0.2cm}
From  the pointwise bound \eqref{ineqmag}, one can also get a diamagnetic inequality for the Dirichlet/Neumann realization in $\Omega$ in the form 
\begin{equation}
|\exp (-t H_{A,V})  (x,y) |\leq  \exp (- t H_{0,V}) (x,y) \,,\, \mbox{ for almost all } x,y \mbox{ in } \Omega ,
\end{equation} for the integral kernel of
the Neumann semigroup (as long as it exists), (see Hundertmark-Simon \cite{HuSi2004}, Remark 1.2 (ii)).

\vspace{0.2cm}
For the trace of $ \exp (-t H_{A,V})$, we can get by a general argument\footnote{If $|A\phi|\leq B|\phi|$ and ${\rm Tr } B^*B < +\infty$, then ${\rm Tr }A^*A\leq {\rm Tr}  B^* B $.} in the proof of  Theorem~15.7  in the book of B. Simon \cite{Si1979}:
 \begin{equation}\label{eq:2.3}
{\rm Tr} \ (\, e^{-t H_{A,V}}) \leq {\rm Tr\,} \ (e^{-t H_{0,V}} )\,.
\end{equation} 
At last, the case of Aharonov-Bohm operator (i.e. when $A$ is an Aharonov-Bohm potential \eqref{defpotAB} ) is  treated in \cite{MOR2005}.

\begin{rem}\label{remark2.1}
	If $\Gamma$ is the complex conjugation, we notice that
	$$
	\Gamma H_{A,V} = H_{-A,V} \Gamma
	$$
	This implies that the eigenvalues of $H_{A,V} $ and $ H_{-A,V}$ coincide (with multiplicity) and consequently
	$$
	{\rm Tr} \ ( e^{-t H_{A,V}}) ={\rm Tr} \ ( e^{-t H_{-A,V}})\,.
	$$
\end{rem}

\vspace{0.2cm}
The existence of complete expansions for the heat kernel is known as $t\rightarrow 0$ at least when $\vec A$ is regular. Some computations rule are given using Boutet de Monvel calculus (see \cite{Bou1971, GiGr1998}). 
The remaining question is then to control the dependence on the magnetic field. If it is clear that the main term is independent of the magnetic field, the question is about the other terms. As already mentioned, a complete magnetic Boutet de Monvel calculus extending \cite{MP2004} could be useful.

\subsection{Magnetic harmonic oscillator in $\mathbb R^2$. }
\subsubsection{Isotropic case.}
In this section, we look at the magnetic harmonic oscillator  in $\R^2$. The Hamiltonian 
$$
H_b:=-\partial_x^2 - (\partial_y + ib x)^2 + x^2 +y^2\,.
$$
defined initially on $\mathcal S(\mathbb R^2)$ has a u nique self-adjoint extension as unbounded operator on $L^2(\mathbb R^2)$.\\
By partial Fourier transform, we first get
$$
\widehat H_b:=-\partial_x^2  -\partial_\eta^2+ (\eta + b x)^2 + x^2\,.
$$
Hence we have now a Schr\"odinger operator without magnetic field but with a new electric potential. We have now to diagonalize the matrix
$$
A_b:=\left(\begin{array}{cc} (1+b^2) &b\\b& 1\end{array}\right)\,.
$$
The eigenvalues are given by
$$
\lambda_\pm = (1+ \frac{b^2}{2}) \pm \frac 12 b \sqrt{4 + b^2}\,.
$$
After a change of variables, we obtain
$$
\check H_b= (-\partial_s ^2 + \lambda_+ s^2) +  (-\partial_z ^2 + \lambda_- z^2)
$$
After a new change of variables, we get
$$
\breve H_b:=\sqrt{\lambda_+} (-\partial_{\tilde s}^2 + \tilde s^2) + \sqrt{\lambda_-} (-\partial_{\tilde z}^2 + \tilde z^2).
$$ 
Now we have
$$
{\rm Tr} \ (e^{- t H_b}) = \Big({\rm Tr\,} e^{-t \sqrt{\lambda_+} (-\partial_{\tilde s}^2 + \tilde s^2)}\Big) \Big({\rm Tr\,} e^{-t \sqrt{\lambda_-} (-\partial_{\tilde z}^2 + \tilde z^2)}\Big).
$$
The computation is easy for the Harmonic oscillator. We obtain:
$$
{\rm Tr\,} e^{-t (-\frac{d^2}{dx^2} + x^2)}= \sum_{n=0}^{+\infty} e^{-(2n+1)t } = \frac{1}{2 \sinh t}\,,\, \forall t>0\,.
$$
So we get
\begin{equation}\label{eq:2.4}
{\rm Tr} \ (e^{- t H_b}) = e^{- (\sqrt{\lambda_+} + \sqrt{\lambda_-}) t} / \Big( (1- e^{-2 \sqrt{\lambda_+}t}) (1- e^{-2 \sqrt{\lambda_-}\,t})\Big)= 1/ (4 \sinh \sqrt{\lambda_+}\,t\,\sinh \sqrt{\lambda_-}\,t)\,.
\end{equation}

\noindent
Coming back to \eqref{eq:2.4}, we now calculate the first terms of the expansion (no logarithm) and get first
$$
{\rm Tr} \ (e^{- t H_b})  = \frac{1}{4t^2}  - \frac {(2+b^2)}{24} + \frac 14 \left( \frac{7}{360} (2+b^2)^2 -\frac{1}{90}\right) t^2 +  \mathcal O (t^4)\,.
$$
Finally, we have 
\begin{equation}
{\rm Tr} \ (e^{- t H_b})  = \frac{1}{4t^2} - \frac {(2+b^2)}{24} + \frac {1}{60} \left(1+  \frac{7}{6} b^2 + \frac{7}{24} b^4\right) t^2 +  \mathcal O (t^4)\,.
\end{equation}
As previously, the coefficient of $t^{-2}$ is independent of $b$. For the relative trace, we get
\begin{equation}
{\rm Tr} \ (e^{- t H_b})  -  {\rm Tr} \ (e^{- t H_0}) = - \frac {b^2}{24} +  \mathcal O (t^2)\,.
\end{equation}

\subsubsection{Anisotropic case}
Now, we look at the same problem in the anisotropic case  in $\R^2$.  We consider the operator
$$
H_b:=-\partial_x^2 - (\partial_y + ib x)^2 + k_1^2 x^2 + k_2^2 y^2\,.
$$
By partial Fourier transform, we first get
$$
\widehat H_b:=-\partial_x^2  - k_2^2 \partial_\eta^2+ (\eta + b x)^2 + k_1^2 x^2\,.
$$
Hence we have now a Schr\"odinger operator without magnetic field but with a new electric potential.  A dilation in $\eta$ leads to
$$
\breve H_b:=-\partial_x^2  -  \partial_\eta^2+ (k_2 \eta + b x)^2 + k_1^2 x^2\,.
$$
We have now to diagonalize the matrix
$$
A_b:=\left(\begin{array}{cc} (k_1^2+b^2) &b k_2 \\b k_2 & k_2^2\end{array}\right)\,.
$$
The eigenvalues are given by
$$
\lambda_\pm = \frac 12 (k_1^2+k_2^2 +b^2) \pm \frac 12\sqrt {\left(\big((k_1-k_2)^2 + b^2\big) \big((k_1+k_2)^2 + b^2\big)   \right)}\,.
$$
As previously, we get
\begin{equation}\label{anisotropic}
	{\rm Tr} \ (e^{- t H_b}) = e^{- (\sqrt{\lambda_+} + \sqrt{\lambda_-}) t} / \Big( (1- e^{-2 \sqrt{\lambda_+}t}) (1- e^{-2 \sqrt{\lambda_-}\,t})\Big)= 1/ (4 \sinh \sqrt{\lambda_+}\,t\,\sinh \sqrt{\lambda_-}\,t)\,.
\end{equation}
As expected the coefficient of $t^{-2}$ is independent of $b$.  Indeed, we  have
$$
\lambda_+ \lambda_- = k_1^2 k_2^2\,,
$$
which is independent of $b$, so we get
$$
{\rm Tr} \ (e^{- t H_b})  = \frac 14 (\lambda_+ \lambda_-)^{-1/2} t^{-2} - \frac{1}{24} (\lambda_+ \lambda_-)^{-1/2}  (\lambda_++\lambda_-) + \mathcal O (t^2) \,.
$$
Finally, we have for the relative trace
$$
{\rm Tr} \ (e^{- t H_b}) - {\rm Tr} \ (e^{- t H_0})= - \frac{1}{24\, k_1k_2} \, b^2 + \mathcal O (t^2).
$$
Up to a sign error we recover Odencrantz formula (cf \cite{Od1988}, Theorem 2).
\begin{rem}
In higher dimension, 
 K. Odencrantz \cite{Od1988} gives the two first terms and considers, keeping the magnetic field constant, more general  electric potentials $V$ (including $V$ quadratic and $V$ quartic). In particular, he gives the main term for 
 ${\rm Tr} \ (e^{- t H_b}) - {\rm Tr} \ (e^{- t H_0})$  as $t \rightarrow 0^+$.  The verification of the corresponding diamagnetic inequality  for the trace of the heat kernel  is easy.
 \end{rem}

\subsection{Aharonov-Bohm harmonic oscillator in $\mathbb R^2$}
 
  Let $A_\nu (x)$ be the Aharonov-Bohm potential in $\R^2$ defined by
 \begin{equation}\label{eq:def-F}
\begin{gathered}
A_\nu (x)= \nu \mathbb F(x),\\
\mathbb F(x)=\begin{pmatrix}-\displaystyle\frac{x_2}{|x|^2}  , \frac{ x_1}{|x|^2}\end{pmatrix}.
\end{gathered}
\end{equation}
In the sense of distribution in $\mathbb R^2$, the magnetic field  is given by $2\pi \nu \ \delta_0$ where $\delta_0$ is the Dirac measure at the origin, (see A. Hansson \cite{Ha2005}).
We now consider the Hamiltonian  in $\mathbb R^2$
$$
H_{\nu,\beta} :=-\Delta_{A_\nu} + \beta r^2\,,
$$
where $-\Delta_{A_\nu}$ is the magnetic Laplacian   and $\beta >0$.
Its spectrum is given by
\begin{equation}\label{eq:2.2}
E(m,n)=2\sqrt{\beta}(1+|m-\nu|+2n),\quad (m\in\Z,~n=0,1,2,\cdots),
\end{equation}
and we have the following formulas 
$$
{\rm Tr\,} e^{-t H_{\nu,\beta}}={\rm Tr\,} e^{- 2 \sqrt{\beta} t H_{\nu,1}}\,,
$$
 Since
\begin{eqnarray}
{\rm Tr\,} e^{-s H_{\nu,1}}&=& e^{-1}\Big( \sum_{m\in \mathbb Z} e^{-s |m-\nu|}\Big)  \Big( \sum_{n\geq 0} e^{-2 s n}\Big)\, \nonumber \\
&=& e^{-1} \frac{\cosh ( s (\{\nu\}-1/2))}{\sinh (s/2)} \frac{1}{1-e^{-2s}}\,,
\end{eqnarray}
where $\{\nu\} = \nu - \lfloor \nu \rfloor$ is the fractional part of $\nu$, (we refer to \eqref{Tdenu} for details), we get immediately:
\begin{equation}
{\rm Tr\,} e^{-t H_{\nu,\beta}}= e^{-1} \ \frac{\cosh ( 2\sqrt{\beta}t (\{\nu\}-1/2))}{\sinh (\sqrt{\beta}t)} \frac{1}{1-e^{-4\sqrt{\beta}t}}\,.
\end{equation}

\subsection{Magnetic Laplacian on the sphere $S^2$}\label{sectionsphere}
 In this situation, the magnetic field is the curvature of a connexion, (see \cite{BeCoCo1998,CoTo1993} ). One can in particular consider the constant magnetic field case  but  its total flux should be an integer. 
For the computation of trace formulas in this case,   \cite{BeCoCo1998} is referring to  \cite{Ku1982}   (formula (7.4)) who sends to \cite{Gi1975}  Gilkey's theorem which  is quite general. We propose now a more explicit expansion.

\vspace{0.2cm}
In the geodesic polar coordinates $(r,\theta)$ with $0<r<\pi\,,\, 0< \theta <2\pi$, the metric is given by $dr^2 +\sin^2 r \, d\theta^2$, the volume form is 
$$\Omega_0:=\sin r \, dr \wedge d\theta$$  and the Hamiltonian reads
$$
-\partial_r^2 - \frac{1}{\sin^2r} \partial_\theta^2 -\cot r\, \partial_r - \frac{2i}{\sin r} f'(r) \partial_\theta + f'(r)^2\,,
$$
where the $1$-form is $\alpha= f'(r) \sin r \,d\theta$.

\vspace{0.2cm}
This leads for the computation of the spectrum to the analysis of the family of the Dirichlet (pour $n\neq 0$)  (Neumann pour $n=0$)) realization  in $ (0,\pi)$ of
$$
-\partial_r^2  -\cot r\, \partial_r + (\frac{n}{\sin r}+ f'(r) )^2 \,.
$$
In the constant magnetic field case, we have  $$ f(r)= \frac m2 (1-\cos r)$$
 with $m\in \mathbb Z$ (see also \cite{KoTa2019}). The spectrum is given by 
 $$
 \nu_{m,j} = j(j+1) + \frac{|m|}{2} (2j+1)\,,\, j=0,1,\cdots $$
 each eigenvalue having multiplicity  $|m|+ 2j+1$.

\vspace{0.5cm}
\par\noindent 
In the following, we assume for simplicity that $m \geq 0$. The trace  of the heat operator $e^{- t H_m}$ is given by:
\begin{eqnarray}\label{Trsphere}
{\rm Tr} \ (e^{- t H_m}) &=&  \sum_{n=0}^{+\infty} (m+2n+1) \ e^{-t\left( n(n+1) + \frac{m}{2}(2n+1)\right)}. \nonumber \\
                        &=&  \left( \sum_{n=0}^{+\infty} (m+2n+1) \ e^{-t\left(n  +\frac{m+1}{2}\right)^2} \right)  \ e^{t \left(\frac{m^2+1}{4}\right)}.
\end{eqnarray}
We can actually get an explicit expression of ${\rm Tr} \ ( e^{- t H_m}) $ in terms of the so-called Jacobi partial theta function \cite{BK2016}, which are defined by:
\begin{align}\label{def_Fdell}
	F_{d,\ell}\left(z;\tau\right):=\sum_{n \geq 0} \zeta^{\ell n+d}q^{(\ell n+d)^2}\,,
\end{align} 
where $d\in\Q^+$, $\ell\in\N$,  $\zeta:=e^{2 \pi i z}$ with $z \in \C$ and $q:=e^{2\pi i \tau}$ with $\tau \in \C^+$, (the complex upper half-plane).  In particular, one gets:
\begin{equation}\label{derivJacobi}
F_{d,\ell}'\left(0;\tau\right) = 2i \pi \sum_{n \geq 0}  (\ell n+d)\ q^{(\ell n+d)^2},
\end{equation}
where $F_{d,\ell}'\left(z;\tau\right)$ stands for the derivative with respect to $z$.
Thus, using (\ref{Trsphere}) and (\ref{derivJacobi}), we easily see that:
\begin{equation}\label{traceJacobi}
 {\rm Tr} \ (e^{- t H_m}) = \frac{1}{i\pi} \ F_{\frac{m+1}{2},1}'\left(0; \frac{it}{2\pi}\right) \ \ e^{t \left(\frac{m^2+1}{4}\right)}.
\end{equation}

To find the asymptotic expansion of ${\rm Tr} \ (e^{- t H_m})$  as $t \to 0^+$, we need the following result which is proved implicitely in (\cite{BK2016}, Theorem 4.1):

\begin{lemma}\label{lemma0.1}
	For $|z|< 1/(4 \ell)$, we have the following asymptotic expansion as $|\tau| \to 0$, 
	$$
	F_{d,\ell}\left(z;\tau\right)=\sum_{j\geq 0}\frac{(2\pi i\ell z)^j}{j!}\left(\frac{\Gamma\left(\frac{j+1}{2}\right)}{2(-2\pi i\ell^2 \tau)^{\frac{j+1}{2}}}
	-\sum_{k=0}^N\frac{\left(2\pi i\ell^2 \tau\right)^k}{k!}\frac{B_{2k+j+1}\left(\frac{d}{\ell}\right)}{2k+j+1}\right)+O\left(|\tau|^{N+1}\right).
	$$
	More precisely, for any compact set $K \subset D(0, \frac{1}{4 \ell})$, the above series converges normally and the remaining term is uniform with respect to $z$. Here $B_n(x)$ denotes the $nth$ Bernoulli polynomial given by 
	$$
	B_n(x) = \sum_{p=0}^n \binom{n}{p} B_p\  x^{n-p},
	$$
where $B_p$ are the Bernoulli numbers.
\end{lemma}

\par\noindent
Now, using  Cauchy's formula: 
$$
F_{\frac{m+1}{2},1}'\left(0; \frac{it}{2\pi}\right) = \frac{1}{2i\pi}\ \int_{\mathcal{C}(0, \frac{1}{8})} \ \frac { F_{\frac{m+1}{2},1}\left(z; \frac{it}{2\pi}\right)}{z^2} \ dz,
$$
and (\ref{traceJacobi}), we get immediately the following asymptotic expansion as $t \to 0^+$:
\begin{equation}
 {\rm Tr} \ (e^{- t H_m}) = \left( \frac{1}{t}  - \sum_{k=0}^N   \frac{(-t)^k}{(k+1)!} \ B_{2k+2} (\frac{m+1}{2})    +O\left(t^{N+1}\right) \right) \ 
 e^{t \left(\frac{m^2+1}{4}\right)}.
\end{equation}
Using Mathematica, we get  :  
\begin{equation}\label{tracesphere}
{\rm Tr} \ (e^{- t H_m}) = \frac{1}{t}+\frac{1}{3}+\left(\frac{1}{15}-\frac{m^2}{24}\right) t+\left(\frac{4}{315}-\frac{m^2}{40}\right) t^2+\frac{\left(128-432 m^2+49
		m^4\right) t^3}{40320}+O(t^4 )
\end{equation}

We emphasize that no logarithmic term appears in the expansion.  In particular, the contribution of the magnetic field appears in the coefficient of $t$. We get for the relative trace 
\begin{equation}
	{\rm Tr} \ (e^{- t H_m})  - {\rm Tr} \ (e^{- t H_0})=  -\frac{m^2}{24} t+ \mathcal O (t^2)\,.
\end{equation}
Here the first term in the right hand side is negative for $t>0$  as predicted by \eqref{eq:2.3}.

\section{On heat expansions for the D-to-N--magnetic operators, general properties}
\subsection{Diamagnetism}

\begin{equation}\label{diamaDN}
	{\rm Tr} \ (e^{-t \Lambda_{A,V}} ) \leq  {\rm Tr} \  ( e^{-t \Lambda_{0,V}}).
\end{equation}	
The first way to get (\ref{diamaDN}) is to use the  general argument given in the proof of Theorem 15.7  in the book of  B. Simon \cite{Si1979}. We can also use the comparaison between the distribution kernels of $\Lambda_{A,0}$ and $\Lambda_{0,0}$, (see  \cite{EO2022}). In this paper, the authors assume that $A$ belongs to $L^2_{loc}$ and that the boundary of $\Omega$ is Lipschitz.

\vspace{0.2cm}
The proof in the case of Aharonov-Bohm  potentials (with poles avoiding the boundary) follows the same line if we use the Friedrichs extension when defining the magnetic D-to--N map
 and use the characterization of the domain as discussed for example in the work of Lena \cite{Le2015}, (but look at two cases depending on an integer (renormalized) flux or not and one should add what results of Hardy's inequality). This  does not change anything at  the boundary. The consequences will also be that the local coefficients will not see the Aharonov-Bohm potentials (more generally magnetic potentials whose corresponding magnetic fields are compactly supported) and this will lead, as $t\rightarrow 0^+$, to
\begin{equation}
	{\rm Tr} \ (e^{-t \Lambda_{A,V}} ) -  {\rm Tr} \  ( e^{-t \Lambda_{0,V}}) = \mathcal O (t)\,.
 \end{equation}
This will be further discussed in the next subsection and explains why we will focus on the computation of the coefficient of $t$ in the expansion of  ${\rm Tr} \ (e^{-t \Lambda_{A,V}} )$  when explicit computations can be done on examples.
 
\begin{rem}
If $\Gamma$ is the complex conjugation, we notice (see Remark \ref{remark2.1})  that
$$
\Gamma \Lambda_{A,V} =\Lambda_{A,V} \Gamma \,.
$$
This implies
$$
 {\rm Tr}  \exp (-t \Lambda_{A,V})= {\rm Tr}  \exp (-t \Lambda_{(-A,V)})
 $$
 and this explains below the parity in $b$ observed below for some of the expansions.
\end{rem}



\subsection{Small $t$ expansions.}
The paper by Liu-Tan  \cite{LiuTan2023} states: 
\begin{thm}\label{thLT}
Suppose that $(\Omega,g)$ is a smooth compact Riemannian manifold of dimension $n$ with smooth boundary $\partial \Omega$. Let $A\in X(\Omega)$ (i.e. a smooth $1$-form on $\Omega$) and $V\in C^\infty(\Omega)$ be the magnetic vector potential and the electric potential. Let $\{\lambda_k\}$ be the eigenvalues of the magnetic Dirichlet-to-Neumann map $\mathcal M$. Then the trace of the heat kernel associated with $\mathcal M$ admits the asymptotics
$$
\sum_j e^{-t \lambda_j} =\sum_{k=0}^{n-1} a_k t^{-n+k+1} + o(1)\mbox{ as } t \to 0^+\,,
$$
where the coefficients are "local" invariants which can be explicitly computed. Moreover the coefficients satisfying $k\leq \inf (n-1,4)$ are independent of the magnetic field.
\end{thm}
When $n=2$, this gives only the information for the main term which depends only on the main term of the D-to-N operator which is independent of the magnetic potential. In Remark~1.5 of \cite{LiuTan2023} it is observed that $a_1=0$.\\ 
For $a_0$, we have the Weyl term relative to the square-root of the Laplace Beltrami   operator on  $\partial \Omega$.  
Note that $o(1)$ can be replaced by $\mathcal \mathcal O(t\log \frac 1t)$ and that this coefficient is "local".
 Hence the first non local  coefficient  where the magnetic potential can play a role is the coefficient of $t$.

\subsection{Extension of Theorem \ref{thLT} to irregular potentials}
In this subsection, we sketch how we can relax the assumptions of regularity on $\Ab$ and $V$ far from $\Omega$ keeping the property
 that $H_{A,V}$ is well defined through the associated variational form.  This is in particular the case when the magnetic potential is a sum of Aharonov-Bohm potentials with singularities in $\Omega$.
 So we assume that the Dirichlet Laplacian is well defined as a self-adjoint operator (using Lax-Milgram) and that the form domain of $H_{A,V}$
  is contained in $H^1(\Omega)$. So for the problem \eqref{Dirichlet}, the application $f \mapsto u$ is continuous from $H^\frac 12(\partial \Omega)$ into $H^1(\Omega)$.

\vspace{0.2cm}\noindent
We now assume that for some $T>0$  $\Ab$ and $V$ are in $C^\infty(\overline{\Omega_T})$, where
  $$
  \Omega_T=\{x\in \Omega\,,\, d(x,\partial \Omega) < T\}\,.
  $$
  With $T$ small enough, we can assume that $\overline{\Omega_T}$ is diffeomorphic to $[0,T] \times \partial \Omega$.
 We have the following lemma: 
 \begin{lemma}\label{reg}
   For $0 < \epsilon < T$, the maps $f \mapsto u_{| \Gamma_\epsilon}$ and $f \mapsto (\nu_\epsilon\cdot \nabla u)_{| \Gamma_\epsilon}$ are continuous from $H^{\frac 12}(\partial \Omega)$ into $H^{k}(\Gamma_\epsilon)$ (for any $k$), where
  $$ \Gamma_\epsilon:= \{x\in \Omega\,,\, d(x,\partial \Omega) =\epsilon \}\,.
  $$
 \end{lemma}

\vspace{0.2cm}
The proof is a direct consequence of the fact that $H_{A,V} \ u=0$ in $\Omega_T$, (and consequently $C^\infty$ there), and the ellipticity property in
$\Omega_T$. \\

We can then follow verbatim the proof given by Lee--Uhlman (\cite{LeeUh1989}, Proposition 1.2) which is a consequence  of a pseudo-differential factorization (\cite{LeeUh1989}, Proposition 1.1). 
Note that we need the condition that $0$ is not in the spectrum of $H_{A,V}$.

\vspace{0.2cm}
Hence, Theorem \ref{thLT} holds under the assumptions of this subsection. For the case with magnetic field and electric field we can refer to the proof of Proposition 3.1 in Liu--Tan (\cite{LiuTan2023}, Prop. 3.1).

\begin{rem} The proof proposed in Lee--Uhlman \cite{LeeUh1989} gives only the statement with a weaker notion of regularizing operator. We say that the operator is weakly regular if it is continuous from $H^\frac  12(\partial \Omega)$  to $H^k(\partial \Omega)$ for any $k\geq 0$.
\end{rem}

\begin{rem} A variant of the above proof is to combine Lemma \ref{reg} with the introduction of  a Dirichlet-to-Neumann operator in  $\Omega_\epsilon$. We are 
now in the regular case. We can take as traces on $ \partial \Omega_\epsilon$, $f_\epsilon = u_{|\Gamma_\epsilon}$ and  $f$ sur $\partial \Omega$.  Actually, we do not need Lemma \ref{reg},  but we can use that for the  D-to-N operator  with two components (which can be considered as a $2\times 2$-system), the off-diagonal terms are strongly regularizing, as it results of the Boutet de Monvel calculus.
It could actually be interesting to develop a magnetic calculus for the problem with boundary (see Mantoiu-Purice \cite{MP2004} with the hope that the contribution of the magnetic field will appear more directly; see also Helffer-Purice \cite{HP2010} for a connected problem).

\end{rem}

\subsection{Coming back to the trace}
In the case $d=2$, with an Aharonov--Bohm potential, and admitting that the existence of the trace expansion is properly proved,  it is interesting to determine where, in the expansion,  the flux created at the singularity appears. The "local" coefficients of the expansion do not see the flux, since except at the singularity which is assumed to be inside the domain, one can gauge away the magnetic potential. Hence, the first non local coefficient would be interesting to analyze. We will come back to this question for particular cases where an explicit computation is possible.

\section{Aharonov--Bohm effect for the trace of the   D-to-N magnetic operator: The case of  the disk.}\label{ABsection}

Note that the Friedrichs extension of the
Aharonov--Bohm Hamiltonian on a disk (that we consider here) and other selfadjoint realizations are discussed by J. F. Brasche and M. Melgaard in  \cite{BM2005}. See also R.~Frank and A.M. Hansson  \cite{FrHa2009}.

\subsection{Framework}
We consider for $\nu \in \mathbb R$  the magnetic Aharonov-Bohm-potential  (in short AB-potential) in $\mathbb R^2\setminus \{0\}$ , (see \cite{AB1959}):
\begin{equation}\label{defpotAB}
A_\nu (x,y) = \frac{\nu}{r^2} (-y,x)\,.
\end{equation}
This magnetic potential creates a flux $2\pi \nu$ around $0$. In the distributional sense we have  in the distributional sense.
$$
{\rm curl}\,A_\nu = \delta_0\,.
$$
We would like to analyze the D-to-N operator denoted
$$
\mathcal M_\nu :=\Lambda_{A_\nu}
$$ 
(in order to simplify the notation) associated with the magnetic Laplacian in the disk $D(0,1)\subset \mathbb R^2$. Then
 the spectrum of $\mathcal M_\nu$ is given by:
 \begin{equation}
 \lambda_k (\nu) =|k-\nu |\mbox{ for } k\in \mathbb Z
 \end{equation}
Except for $\nu\equiv \frac 12$ modulo $\mathbb Z$, each eigenvalue has multiplicity $1$. 
 We are interested in
 \begin{equation} \label{traceAB}
 (t,\nu) \mapsto {\rm Tr } (e^{-t \mathcal M_\nu} )= \sum_{k\in \mathbb Z} e^{-t |k-\nu|}\,.
 \end{equation}
This function is $1$-periodic with respect to $\nu$ and from now on we assume that $\nu \in [0,1)$.\\

\subsection{Explicit computations}
By elementary computation, we get, assuming that $\nu \in [0,1[$,

$$
{\rm Tr } (e^{-t \mathcal M_\nu} )=\sum_{k=0}^{+\infty} e^{-t (k+\nu)} + \sum_{k=0}^{+\infty} e^{-t (k-\nu)} - e^{t\nu}
$$
i.e.
\begin{equation}\label{Tdenu}
\mathfrak T(\nu):= {\rm Tr } (e^{-t \mathcal M_\nu} )= \frac{2 \cosh t\nu}{1-e^{-t}} -e^{t\nu} =  \frac{\cosh ( t (\nu-1/2))}{\sinh (t/2)}\,.
\end{equation}
From this formula, we see immediately that $\nu \mapsto \mathfrak T(\nu)$ is symmetric with respect to $\frac 12$ and attains its minimum at $\nu=\frac 12$.
We also have the asymptotic
\begin{equation}\label{Tdenua}
\mathfrak T(\nu) = \frac 2t + (\frac 16 -\nu + \nu^2) t + \mathcal O (t^3)\,.
\end{equation}
We note that the first term where the magnetic flux appears is $(\frac 16 -\nu + \nu^2) t$. \\
We can also write:
$$
\mathfrak T(0) -\mathfrak T(\nu) =2 (\cosh (t/2) - \cosh ( t (\nu-1/2)) /\sinh (t/2) = 4 \sinh (t\nu) \sinh ((1-\nu)t)/ \sinh (t/2) \geq 0\,.
$$
This gives in particular the diamagnetic property which holds in a more general situation (see the papers by Ter Elst--Ouhabaz \cite{EO2022} and references therein for the regular case).

\subsection{Generalization by changing the metrics}
One can generalize the preceding computations by working in the ball but by changing the metrics (in the disk direction). For this we introduce a function $\theta (r)$ such that
$$
\theta >0 \mbox{ on } (0,1) \mbox{ and } \theta'(0) =1\,.
$$
In this case, the Laplacian (without magnetic field) reads in polar coordinates $(r,t)$
$$
-\Delta = - \frac{\partial^2}{\partial r^2} - \frac{\theta'(r)}{\theta(r)} \,\frac{\partial }{\partial r} - \frac{1}{\theta(r)^2} \frac{\partial^2}{\partial t^2}\,.
$$

We then keep the same magnetic Aharonov--Bohm $1$-form potential $\nu dt$ as above and (as computed in \cite{CPS2022}, Subsection 4.1) we introduce $-\Delta_{A_\nu}$ 
$$
-\Delta_{A_\nu}=- \frac{\partial^2}{\partial r^2} - \frac{\theta'(r)}{\theta(r)} \,\frac{\partial }{\partial r} - \frac{1}{\theta(r)^2} \big(\frac{\partial}{\partial t} - i \nu\big) ^2\,.
$$
Let us compute the Steklov eigenvalues. We have to solve in $(0,1)$ (with $u$ bounded near $0$)
\begin{equation}\label{eq:AB1}
u'' + \frac{\theta'}{\theta} u' - \frac{(k-\nu)^2}{2} u=0 \mbox{ and } \, u'(1)= \sigma u(1)\,.
\end{equation}
We solve \eqref{eq:AB1} without the Robin condition and get
$$
u_k (r) = e^{|k-\nu| \int_1^r \frac{1}{\theta(s)} ds}\,.
$$
We then recover the  Steklov's eigenvalue $\sigma_k$ by writing
\begin{equation}
\sigma_k = \frac{u'_k(1)}{u_k(1)} = |k-\nu|/\theta(1)\,.
\end{equation}
The computations are unchanged as for Aharonov-Bohm case after the change of variable $t \mapsto t/\theta(1)$.  Denoting by $\mathcal{M}_{\nu}^{\theta}$ the associated D-to-N map, we get: 

\begin{equation}\label{Tdenutheta}
	{\rm Tr } (e^{- t \mathcal M_\nu^{\theta}} )= \frac{\cosh ( \frac{t}{\theta(1)} (\nu-1/2))}{\sinh (\frac{t}{2\theta(1)})}\,.
\end{equation}

\subsection{Generalization to  the cylinder}\label{Cylindersection}

We consider the cylinder $M=(-1,1)\times S^1$ and let
\begin{equation}
H_\nu  = -\frac{\partial^2}{\partial x^2} + (i^{-1} \frac
{\partial}{\partial \theta} -\nu)^2 \mbox{ with } \nu\in (0,1) ,
\end{equation}
be the magnetic Laplacian. It is proved in \cite{PrSa2023} that the  Steklov spectrum  of the D-to-N map $\Lambda_\nu$ associated with $H_\nu$ is given by:
\begin{equation}
\begin{array}{ll}\label{vp}
\lambda_k^+(\nu )&= (k-\nu ) \coth (k-\nu ) \ , \\
\lambda_k^-(\nu )&= (k-\nu ) \tanh (k-\nu ) \ ,\ k \in \mathbb Z.
\end{array}
\end{equation}
Let us compute  the asymptotics as $t\rightarrow 0^+$ of

\begin{equation}
	{\rm Tr} (e^{-t \Lambda_\nu})= \sum_{k \in \Z} e^{-t \lambda_k^+(\nu)} +  \sum_{k \in \Z} e^{-t \lambda_k^-(\nu)} := S^+(t,\nu) + S^-(t,\nu)\,.
\end{equation}
For instance, we get for $S^-(t,\nu)$,
\begin{equation}
S^-(t,\nu) = \sum_{k \in \Z} e^{-t |k-\nu|} +   \sum_{k \in \Z} ( e^{-t \lambda_k^-(\nu)} - e^{-t |k-\nu|}),
\end{equation}
and using \eqref{traceAB} we obtain:
\begin{equation} \label{decomposition}
S^-(t,\nu) = {\rm Tr} (e^{-t \mathcal M (\nu)}) + R^-(t,\nu)\,,
\end{equation}
where 
\begin{equation}
R^-(t,\nu)  =  \sum_{k \in \Z} \  \sum_{p \geq 1}  \frac{(-t)^p}{p!} \left( (\lambda_k^{-}(\nu))^p - |k-\nu|^p \right).
\end{equation}
A straightforward calculation shows there exists a positive constant $C$ such that
\begin{equation}
	|  (\lambda_k^{-}(\nu))^p - |k-\nu|^p | \leq C \ p \ (|k|+1)^p e^{-2  |k|}\ , \ k \in \Z.
\end{equation}
We deduce that 
\begin{eqnarray*}
	\sum_{p \geq 1}   \frac{t^p}{p!} \ \sum_{k \in \Z} \  \left| (\lambda_k^{-}(\nu))^p - |k-\nu|^p \right| &\leq& C \ 
	\sum_{p \geq 1}   \frac{t^p}{(p-1)!} \sum_{k \in \Z} \  (|k|+1)^p e^{-2 |k|} \\
	&\leq& C \ \sum_{p \geq 1}   \frac{t^p}{(p-1)!} \ \frac{p!}{2^p}  < \infty,
\end{eqnarray*}
for $t \in ]0,2[$. We deduce that for such $t$,
\begin{equation} \label{restenegatif}
	R^-(t,\nu)  =   \sum_{p \geq 1} a_p^-  \ t^p ,
\end{equation}
where 
\begin{equation}
	a_p^- =  \frac{(-1)^p}{p!}  \sum_{k \in \Z} \ \left( (\lambda_k^{-}(\nu))^p - |k-\nu|^p \right)
\end{equation}
Then, using (\ref{decomposition}) and (\ref{restenegatif}), one has for $t \in ]0, 2[$,
\begin{equation}
S^- (t, \nu) = {\rm Tr} (e^{-t \mathcal M (\nu)}) +  \sum_{p \geq 1} a_p^-  \ t^p .
\end{equation}
In the same way, one has:
\begin{equation}
	S^+ (t, \nu) = {\rm Tr} (e^{-t \mathcal M (\nu)}) +  \sum_{p \geq 1} a_p^+  \ t^p ,
\end{equation}
where
\begin{equation}
a_p^+ =  \frac{(-1)^p}{p!}  \sum_{k \in \Z} \ \left( (\lambda_k^{+}(\nu))^p - |k-\nu|^p \right)	.
\end{equation}
It follows immediately that for $t \in ]0,2[$,
\begin{equation}
{\rm Tr} (e^{-t \Lambda_\nu}) =  2\  {\rm Tr} (e^{-t \mathcal M (\nu)})  + 
\sum_{p \geq 1} (a_p^+ + a_p^- )  \ t^p .
\end{equation}
Then, thanks to \eqref{Tdenu},  we get a complete expansion of ${\rm Tr} (e^{-t \Lambda_\nu})$ in powers of $t$ and in particular, one gets using (\ref{vp}):
\begin{equation}\label{eq:trdisk}
{\rm Tr} (e^{-t \Lambda_\nu}) =  \frac {4}{t} + 2 \left(\frac 16-\nu + \nu^2 - \breve a_1(\nu) \right) t  +\mathcal \mathcal O(t^2),
\end{equation}
with 
\begin{equation}
\breve a_1(\nu)= \sum_{k\in \Z} |k-\nu| \left( \coth (2 |k-\nu|)-1\right)\,.
\end{equation}

\section{Aharonov--Bohm effect for the trace of the D-to-N magnetic operator-- The case of the annulus}

\subsection{Trace formulas for the full D-to-N map.}
 
Let us consider an annulus $\Omega \subset \R^2$ centerd at this origin with $R_0 =1$ as external radius and $R<1$ as internal radius. Hence, in radial coordinates
$$
\Omega=\{ (r,\theta) | R < r <1\}\,.
$$
The boundary $\partial \Omega$ has two components
$$
\partial \Omega = S^R \cup S^1\,.
$$
Hence the D-to-N-- map reduced to ${\rm Vect}( e^{im \theta})$, ($m\in \mathbb Z$), will be a $2\times 2$ matrix.
We first solve the magnetic Dirichlet problem (with $A=A_\nu$ as in Section \ref{ABsection}  and $H = (D-A_\nu)^2$):
\begin{equation}
	H v=0 \mbox{ in } \Omega\mbox{ and } v =\Psi \mbox{ on } \partial \Omega\,,
\end{equation} 
where
$$
\Psi = \left(\begin{array}{c} \psi^1\\\psi^R\end{array}\right) \in H^{1/2}(S^1) \oplus H^{1/2}(S^R)\,.
$$
Here we follow \cite{PrSa2023} which refers to \cite{FrSc2011}. We proceed by separation of variables. More precisely we observe 
that the D-to-N $\mathcal M$ map commutes with the rotation. Hence we can consider the joint spectrum of $\partial_{\theta}$ and $\mathcal M$ and write
\begin{equation}
v(r,t)= \sum_{m\in \Z} v_m (r) e^{im \theta}\,.
\end{equation}
At the boundary, writing
\begin{equation}
\psi^{1,R}(t) =\sum_{m\in \Z} \;\psi_m^{1,R}\  e^{im \theta} 
\end{equation}
we get for $m \in \Z$,
\begin{subequations}
	\begin{equation}
		- v''_m -\frac 1r v_m' +\frac{1}{r^2} (m-\nu)^2 v_m=0\,,
	\end{equation}
with the boundary conditions
	\begin{equation}
		v_m(1)= \psi_m^1\mbox{ and } v_m(R) =\psi_m^R\,.
	\end{equation}
\end{subequations}
This leads to
	\begin{equation}
v_m(r) = A_m r^{|m-\nu|} + B_m r^{-|m-\nu|},
	\end{equation}
with
	\begin{equation}
		\left\{
\begin{array}{rl}
	A_m + B_m & = \psi_m^1 \\
	A_m R^{|m-\nu|} + B_m R^{-|m-\nu|} &=\psi_m^2\,.
\end{array}
\right.
\end{equation}
We get:
\begin{eqnarray*}
	A_m &=& - \frac{R^{-\m}}{R^{\m} - R^{-\m}} \ \psi_m^1 + \frac{1}{R^{\m} - R^{-\m}} \  \psi_m^R, \\
	B_m &=&  \frac{R^{\m}}{R^{\m} - R^{-\m}} \  \psi_m^1 - \frac{1}{R^{\m} - R^{-\m}} \  \psi_m^R .
\end{eqnarray*}
The D-to-N reduced to ${\rm Vect}( e^{im \theta})$ is given by
\begin{equation}
	\mathcal{M}_m \ 
	\begin{pmatrix}
		\psi_m^1\\[3mm]
		\psi_m^R \\[3mm]
	\end{pmatrix}
=
\begin{pmatrix}
	v_m'(1)\\[3mm]
	- v_m'(R)\\[3mm]
\end{pmatrix}
\end{equation}
Then, we easily see  that the reduced D-to-N map $\mathcal{M}_m$ is a $2\times 2$ matrix given by 
\begin{equation}\label{matrix}
	\mathcal{M}_m \  = \frac{\m}{R^{\m} - R^{-\m}}
\begin{pmatrix}
	- (R^{\m}+R^{-\m}) & 2 \\[5mm]
	\frac{2}{R} & - \frac{1}{R} (R^{\m}+R^{-\m})
\end{pmatrix}
\end{equation}
Then, a straightforward calculation shows us that the two eigenvalues are given by

\begin{eqnarray}
	\lambda_m^{\pm}& =& \frac{|m-\nu|}{2 (R^{|m-\nu|}- R ^{-|m-\nu|})}\left(  - (1+ \frac 1R)( R^{|m-\nu|}+R ^{-|m-\nu|}) \right. \nonumber \\
	& & \left. \mp \sqrt{ (1-\frac 1R)^2 (R^{|m-\nu|}+ R ^{-|m-\nu|})^2 + \frac{16}{R}}\right).
\end{eqnarray}
Now, writing $R=e^{-a}$ with $a>0$, one easily gets the following asymptotics when  $ m \to \pm \infty$,
\begin{eqnarray}
\lambda_m^+ &=& \frac{\m}{R} \  (1+\mathcal O(e^{-2 a|m|})) , \label{asymptotic-a} \\
\lambda_m^- &=& \m \  (1 +\mathcal O(e^{-2 a|m|})). \label{asymptotic-b}
\end{eqnarray}
Now, we can follow exactly the same approach as in Section \ref{Cylindersection}. Assuming that $\nu \in [0,1[$, we see that the trace of $e^{-t \mathcal M_\nu(R)}$ has a complete asymptotic expansion  for $t$ small enough,
\begin{equation}
{\rm Tr} (e^{-t \mathcal M_\nu(R)}) = 	\frac{\cosh ( t (\nu-1/2))}{\sinh \frac{t}{2}} + \frac{\cosh ( \frac{t}{R} (\nu-1/2))}{\sinh \frac{t}{2R}} + \sum_{p \geq 1} a_p(R,\nu) \ t^p,
\end{equation}
with
\begin{equation}
	a_1(R,\nu): = - (1+\frac{1}{R})  \ \sum_{m \in \Z} \m \left( \coth a \m -1 \right),
\end{equation}
where we have set  $$ R= e^{-a} \mbox{ for }a>0.$$
 In particular, one has 
\begin{equation}
{\rm Tr} (e^{-t \mathcal M_\nu(R)}) = \frac{2+2R}{t} + \left( \frac{1+R-6\nu -6R\nu +6\nu^2 +6R\nu^2}{6R} +a_1(R,\nu) \right) t +\mathcal O(t^2).
\end{equation}

\begin{rem}
 Using (\ref{asymptotic-a}) and (\ref{asymptotic-b}), we easily see that, for all $t >0$, 
 \begin{equation}
 {\rm Tr} (e^{-t \mathcal M_\nu(R)}) \to {\rm Tr} (e^{-t \mathcal M_\nu}) \ , \ R \to 0,
 \end{equation}
 where $\mathcal M_\nu$ is the Aharonov-Bohm D-to-N map defined in Section \ref{ABsection}.
 \end{rem}

\subsection{Trace formulas for the partial D-to-N map.}

We are also interested with trace formulas for the \emph{partial} D-to-N maps defined as follows. Let $\Gamma_D$ and $\Gamma_N$ be two open subsets of $\partial \Omega$. We define the partial D-to-N map $\mathcal{M}_{\nu,\Gamma_D,\Gamma_N} $ as the restriction of the global D-to-N map $\mathcal{M}_\nu$ to Dirichlet data given on $\Gamma_D$ and Neumann data measured on $\Gamma_N$. Precisely, consider the Dirichlet problem
\begin{equation} \label{Eq0}
	\left\{ \begin{array}{cc} H u = 0, & \textrm{on} \ \Omega, \\ u = \psi, & \textrm{on} \ \Gamma_D, \\ u = 0, & \textrm{on} \ \partial \Omega \setminus \Gamma_D. \end{array} \right.
\end{equation}
We define $\mathcal{M}_{\nu,\Gamma_D,\Gamma_N}$ as the operator acting on the functions $\psi \in H^{1/2}(\partial \Omega)$ with $\textrm{supp}\,\psi \subset \Gamma_D$ by
\begin{equation} \label{Partial-DNmap}
	\mathcal{M}_{\nu,\Gamma_D,\Gamma_N} (\psi) = \left( \partial_\nu u \right)_{|\Gamma_N},
\end{equation}
where $u$ is the unique solution of (\ref{Eq0}). \\
In particular, if we take $ \Gamma_D = \Gamma_N = S^1$, and if denote 
$\mathcal{M}_{m, 1}$ this partial D-to-N map reduced to $Vect(e^{im\theta})$, it is easy to see that $\mathcal{M}_{m, 1}$  is only the multiplication operator by the first entry of the matrix $\mathcal{M}_m$ given in (\ref{matrix}), i.e we have:
\begin{equation}
	\mathcal{M}_{m, 1} \ e^{im\theta} = - \frac{\m}{R^{\m} - R^{-\m}} \ (R^{\m}+R^{-\m}) \ e^{im\theta}
\end{equation}
As previously and assuming that $\nu \in [0,1[$, it follows that the heat trace for the partial D-to-N map  $\mathcal{M}_{\nu,S^1,S^1}$ has a complete asymptotic expansion  for $t$ small enough,
\begin{equation}\label{partial}
	{\rm Tr} (e^{-t  \mathcal{M}_{\nu,S^1,S^1}  (R)   }) = 	\frac{\cosh ( t (\nu-1/2))}{\sinh \frac{t}{2}} + \sum_{p \geq 1} \tilde a_p(R,\nu) \ t^p,
\end{equation}
with
\begin{equation}
	\tilde a_1 (R,\nu) = -  \ \sum_{m \in \Z} \m \left( \coth a \m -1 \right),
\end{equation}
where we have set  $ R= e^{-a}$, $a>0$. In particular, one has 
\begin{equation}
	{\rm Tr} (e^{-t  \mathcal{M}_{\nu,S^1,S^1}(R)}     ) = \frac{2}{t} + \left( \frac{1}{6}  - \nu + \nu^2 +\tilde a_1(R,\nu) \right) t +\mathcal O(t^2).
\end{equation}

\section{ Constant magnetic field in the disk} \label{champconstant}

\subsection{Framework}
In this section, we consider the following magnetic $1$-form defined in the unit disk $D(0,1) \subset \R^2$ by :
\begin{equation}\label{constantb}
	A(x,y) = b \, (-y dx + xdy),
\end{equation}
where $b$ is a fixed constant. The $2$-form $dA= 2b\  dxdy$ is a constant magnetic field of strength $2b$. Note that the magnetic potential $A$ satisfies the Coulomb gauge : ${\rm div} A=0$. 
The magnetic Laplacian associated with  this potential $A$ is given by 
\begin{equation}\label{hamiltonian}
	\Delta_A = (D-A)^2 .
\end{equation}
First, we have to solve:
\begin{equation}
\left\{
\begin{array}{ll}
	\Delta_A \ v =0  & \mbox{in} \  D(0,1) , \\
	v= \Psi  & \mbox{on} \ S^1 .
\end{array}
\right.
\end{equation}
Then, in polar coordinates $(r, \theta)$, the D-to-N map is defined (in a weak sense) by :
\begin{equation}\label{defD-t-N--}
\begin{array}{lll}
\Lambda(b) : H^{\half} (S^1) & \to &    H^{-\half} (S^1) \\
	 \hspace{1.5cm} \Psi  &\to& \partial_r v (r, \theta)|_{r=1} .
\end{array}
\end{equation}

\noindent
Writing 
\begin{equation}
v(r, \theta) = \sum_{n \in \Z} v_n (r) e^{in \theta}\ \ , \ \  \Psi(\theta) = \sum_{n \in \Z} \Psi_n e^{in \theta}, 
\end{equation}
we see, (\cite{CLPS2023}, Appendix B), that $v_n (r)$ solves:
\begin{equation}\label{polarequations}
	\left\{
	\begin{array}{ll}
		- v_n'' (r) - \frac{v_n'(r)}{r} + (br-\frac{n}{r} )^2 v_n (r)= 0   & \mbox{for} \  r \in (0,1) , \\
		v_n (1) = \Psi_n .
	\end{array}
	\right.
\end{equation}
A bounded solution to the differential equation (\ref{polarequations}) is given by, (see \cite{CLPS2023}, Eq. (B.2)):
\begin{equation}\label{solutionpos}
	v_n(r) = e^{-\frac{br^2}{2}} r^n L_{-\half}^n (br^2)   \ \ \mbox{for} \ n \geq 0 ,
\end{equation}
where $ L_{\nu}^\alpha (z)$  denotes the generalized Laguerre function. For $n \leq -1$, thanks to symmetries in  (\ref{polarequations}), we get a similar expression for $v_n(r)$  
changing the parameters $(n,b)$ into $(-n, -b)$.

\vspace{0.3cm}\noindent
We recall, (\cite{MOS1966}, p. 336), that the  generalized Laguerre functions $L_{\nu}^\alpha (z)$ satisfy the differential equation:
\begin{equation}
	z \frac{d^2 w}{dz^2} + (1+\alpha-z) \frac{dw}{dz} + \nu w =  0 ,
\end{equation}
and are given by
\begin{equation} \label{lienM}
	L_{\nu}^\alpha (z) = \frac{\Gamma(\alpha+ \nu+1)}{\Gamma (\alpha+1)\Gamma(\nu+1) }  M(-\nu,\alpha+1,z),
\end{equation}
where $M(a,c,z)$ is the Kummer's confluent hypergeometric functio, (also denoted by $_1F_1 (a,c,z)$ in the literature), defined as 
\begin{equation}
	M(a,c,z) = \sum_{n=0}^{+\infty} \frac{(a)_n}{(c)_n} \ \frac{z^n}{n!}.
\end{equation}
Here $(a)_n = \frac{\Gamma(a+n)}{\Gamma(a)}$ is the so-called Pochhammer's symbol and $c \notin \Z^-$, (see \cite{MOS1966}, p. 262). \\

\noindent
For $0 < a < c$, we have the  following formula (\cite{MOS1966}, p. 274):
\begin{equation}\label{eq:kummer}
M (a,c,z)= \frac{\Gamma(c)}{\Gamma(c-a)\Gamma(a)} \int_0^1 e^{zt} t^{a-1} (1-t)^{c-a-1} dt\,.
\end{equation}
When $z=0$, we remark that the previous integral is the usual  Beta function $\beta(x,y)$ computed   at $(a,c-a)$. It follows that $M(a,c,0)=1$. The derivative of the Kummer's function 
with respect to $z$ is given by:
 \begin{equation}\label{derivM}
 \partial_z M(a,c,z) =\frac{a}{c} \  M(a+1,c+1,z)\,.
 \end{equation}
 More generally, for all $N \in \N$, we get:
 \begin{equation}
 \partial_z^N M(a,c,z) =\frac {(a)_N}{(c)_N}  M(a+N,c+N,z)\,,
 \end{equation}
 (see \cite{MOS1966}, p. 264). Now, writing the Taylor expansion of $M(a,c,z)$ at $z_0=0$ with integral remainder, we easily get:
 \begin{equation}
\Big|M(a,c,z) -  \sum_{r=0}^{N} \frac{(a)_r}{(c)_r} \ \frac{z^r}{r!} \Big| \leq \frac{1}{(N+1)!}\,\frac {(a)_{N+1}}{(c)_{N+1}}   |z|^{N+1} \sup_{s\in [0,z]}| M(a+N+1,c+N+1,s)|.
 \end{equation}
 We observe that
 \begin{eqnarray*}
  \sup_{ s\in [0,z]}| M(a+N+1,c+N+1,s)| & \leq & M(a+N+1,c+N+1,z) \\
  & \leq &  e^{z_+}  M(a+N+1,c+N+1,0) = e^{z_+},
 \end{eqnarray*}
 where $z_+ = {\rm max } (z,0)$.\\
  Hence we finally get, for $0 <a < c$, and $z\in \mathbb R$, 
 \begin{equation}\label{asymptMimp}
  \Big|M(a,c,z) -  \sum_{r=0}^{N} \frac{(a)_r}{(c)_r} \ \frac{z^r}{r!} \Big| \leq \frac{1}{(N+1)!}\,\frac {(a)_{N+1}}{(c)_{N+1}}   |z|^{N+1} e^{z_+} \,.
 \end{equation}

\vspace{0.3cm}\noindent
Now, let us return to the study of the Steklov eigenvalues. Obviously, they are given by
\begin{equation}
	\lambda_n = \frac {v_n' (1)}{v_n (1)} \ \ \mbox{for} \ n \in \Z .
\end{equation}
Thus, using (\ref{solutionpos}) and (\ref{lienM}), we see that the Steklov spectrum is the set:
\begin{equation}\label{spectrumSpos}
	\sigma(\Lambda(b)) = \{\lambda_0(b)\} \cup \{\lambda_n(b) , \lambda_n(-b) \ \}_{n \in \N^*},
\end{equation}
where for $n \geq 0$,
\begin{equation}\label{explicitvp}
\lambda_n (b) =  n - b +  2b \ \frac{\partial_z M(\half, n+1, b)}{M(\half, n+1,b)}.
\end{equation}

\begin{rem}
Note that the case $n=0$  is  unambiguous. Indeed, recalling that the hypergeometric function  $M(\half, 1, b)=  e^{\frac{b}{2}}\ I_0(\frac{b}{2})$ where 
	\begin{equation}\label{BesselI0}
		I_0(z) =\sum_{k=0}^{+\infty}  \frac{z^{2k}}{2^{2k} (k!)^2}
	\end{equation} 
is the modified Bessel function of the first kind of order $0$, we get immediately:
\begin{equation}\label{lambda0}
\lambda_0 (b) = b \ \frac{I_0'(\frac{b}{2})}{I_0(\frac{b}{2})},
\end{equation}
and this last quantity is even with respect to $b$.
\end{rem}


\subsection{Asymptotic expansion  of the relative Steklov heat trace operator.}
In this subsection, we compute  the asymptotics when $t \to 0^+$ of the relative Steklov heat trace operator which is given by:
\begin{eqnarray}
 {\rm Tr\,}(e^{-t\Lambda(b)} - e^{-t\Lambda(0)}) &=& \sum_{n \in \Z} (e^{-t \lambda_n (b)} - e^{-t |n|}) \nonumber \\
                                                 &=&  \sum_{n \in \Z} e^{-t \lambda_n (b)} - \coth (\frac{t}{2}).
                                                 \end{eqnarray}
We assume that the strength of the magnetic field $b \in [-L,L]$, where $L$ is a fixed positive real. We emphasize that in the following, we control uniformly all the remainders appearing in  
the next asymptotic expansions with respect to small $b$. 
\vspace{0.2cm}
\par\noindent
Using (\ref{derivM}), (\ref{asymptMimp}) and (\ref{explicitvp}), a tedious calculus shows that  
\begin{equation}\label{asymppositif}
	\lambda_n (b) = n - b + 2b \left( \frac{1}{2n} + \frac{b-1}{2n^2} + \frac{2b^2- 7b+2}{4n^3} \right) + \mathcal O(\frac{|b|}{n^4}) \ \ , \ \ n \to + \infty .
\end{equation}

\vspace{0.3cm}\noindent
For $n \in \N^*$, we set 
\begin{equation}\label{mun}
	\mu_n (b)= n - b + 2b \left( \frac{1}{2n} + \frac{b-1}{2n^2} + \frac{2b^2- 7b+2}{4n^3} \right) ,
\end{equation}
so that 
\begin{equation} 
 \lambda_n(b)- \mu_n(b) = \mathcal O(\frac{|b|}{n^4}) \ , \  n \to + \infty.
\end{equation}

\vspace{0.3cm}\noindent
Now we can compute  the asymptotics of ${\rm Tr\,}(e^{-t\Lambda(b)} - e^{-t\Lambda(0)}  )$ when $t \to 0^+$. We write 
\begin{equation}\label{trace}
\tr (e^{-t\Lambda(b)}-e^{-t\Lambda(0)}) =  (e^{-t\lambda_0 (b)}-1) + \sum_{n=1}^{+\infty} (e^{-t\lambda_n (b)} - e^{-tn})  + \sum_{n=1}^{+\infty} (e^{-t\lambda_n (-b)}- e^{-tn}).
\end{equation}
First, let us study the first term of the (RHS)  of (\ref{trace}). Using (\ref{BesselI0}) and  (\ref{lambda0}), we get the following asymptotics when $b \to 0$:
\begin{equation}
	\lambda_0(b)= \frac{b^2}{4} + \mathcal O(b^4),
\end{equation}
It follows that the first term of the (RHS)  of (\ref{trace}) satisfies
\begin{equation}\label{trace0}
e^{-t\lambda_0(b)} = 1 -t \lambda_0(b) +\frac{t^2}{2} \lambda_0^2(b) + \mathcal O(b^6 t^3) \ ,\ t \to 0^+.
\end{equation}
Now, let us study  the second term of the (RHS)  of (\ref{trace}). To simplify the notation in the beginning of this calculus, we will  write $\lambda_n = \lambda_n (b)$ and $\mu_n = \mu_n (b)$.

\par\noindent
One has:
\begin{eqnarray*}
		\sum_{n=1}^{+\infty} (e^{-t\lambda_n(b)} - e^{-tn}) &=& \sum_{n=1}^{+\infty} (e^{-t\mu_n} \ e^{-t(\lambda_n - \mu_n)} - e^{-tn})\\ 
                           	&=& \sum_{n=1}^{+\infty} \left[e^{-t\mu_n} \left( 1-t(\lambda_n - \mu_n) + \frac{t^2}{2} (\lambda_n - \mu_n)^2 +\mathcal O(\frac{b^3t^3}{n^{12}}) \right) - e^{-tn} \right]   \\   
                           	&=& \sum_{n=1}^{+\infty} (e^{-t\mu_n}- e^{-tn})  -t \sum_{n=1}^{+\infty} e^{-t\mu_n} (\lambda_n - \mu_n) 
                           	+ \frac{t^2}{2} \sum_{n=1}^{+\infty} e^{-t\mu_n} (\lambda_n - \mu_n)^2 + \mathcal O(|b|^3 t^3)  \\
                           	&:=& (1)+(2)+(3) + \mathcal O(|b|^3 t^3).
\end{eqnarray*}

\vspace{0.3cm}\noindent
First, let us give the asymptotic expansion of $(1)$. We set 
\begin{equation}
 \gamma_n = \mu_n -n+b = 2b \left( \frac{1}{2n} + \frac{b-1}{2n^2} + \frac{2b^2-7b+2}{4n^2} \right)
\end{equation}
One has:
\begin{eqnarray*}
	(1) &=& \sum_{n=1}^{+\infty} (e^{-t(n-b +\gamma_n) } - e^{-tn}) \\
	&=& e^{tb} \sum_{n=1}^{+\infty} (e^{-t(n +\gamma_n) } - e^{-tn}) + (e^{tb} -1) \sum_{n=1}^{+\infty} e^{-tn} \\
	&=& e^{tb}\sum_{n=1}^{+\infty} e^{-tn} ( e^{-t\gamma_n } - 1) + \frac{ e^{tb} -1}{e^t -1} \\
	&=&  e^{tb} \sum_{n=1}^{+\infty} e^{-tn} \left( - t\gamma_n + \frac{t^2}{2} \gamma_n^2 + \mathcal O(\frac{|b|^3 t^3}{n^3}) \right) + \frac{ e^{tb} -1}{e^t -1}    \\
	&=&  e^{tb} \sum_{n=1}^{+\infty} e^{-tn} \left( - t\gamma_n + \frac{t^2}{2} \gamma_n^2 \right)  + \frac{ e^{tb} -1}{e^t -1} + \mathcal O(|b|^3 t^3).
\end{eqnarray*}
Using Mathematica, we obtain the following asymptotics when $t \to 0^+$:
\begin{equation}\label{Asympt1}
	(1) = b + bt \log t +C(b)\ t + b t^2 \log t + D (b)\ t^2 + \frac{1}{4} (2b-b^2) t^3 \log t + \mathcal O(|b| t^3),
\end{equation}
where 
\begin{equation*}
C(b) = \frac{1}{6} \bigg( -3 b+3 b^2+ ( b - b^2) \pi ^2- (6b -21 b^2 +6 b^3 )\ \zeta(3) \bigg) ,
\end{equation*}
and 
\begin{eqnarray*}
	D(b) &=& \frac{1}{7560}\bigg( -10710 \,b+5670\, b^2+1260 \,b^3+ ( 1260 \, b -2520\, b^2) \pi ^2 + (126\, b^2 -378 \, b^3 +126 \, b^4 ) \pi ^4   \\ 
		& & + ( 4 \, b^2 -28 \, b^3 +57\, b^4 -28 \,b^5 +4 \,b^6 ) \pi ^6 - (15120 \,b^2 -34020\, b^3 +7560 \,b^4 )\ \zeta(3)  \\
        & &  - (7560 \,b^2 -34020\, b^3  +34020 \, b^4 -7560\, b^5 ) \  \zeta(5) \bigg).
\end{eqnarray*}

\vspace{0.3cm}\noindent
Secondly, let us specify the asymptotics of $(2)$. One has:
\begin{eqnarray*}
(2)    &=& 	-t \sum_{n=1}^{+\infty} e^{-t\mu_n} (\lambda_n - \mu_n) \\
       &=& -t \sum_{n=1}^{+\infty} \left( 1 - t\mu_n +R_n(t) \right) (\lambda_n - \mu_n),
\end{eqnarray*}
where
\begin{equation}
	R_n(t) = (t\mu_n)^2 \int_0^1 (1-s) e^{-st\mu_n} \ ds
\end{equation}
and satisfies 
\begin{equation}
	|R_n (t) | \leq \frac{t^2 \mu_n^2}{2}.
\end{equation}
It follows that 
\begin{equation} \label{Asympt2}
	(2)= -t \sum_{n=1}^{+\infty} (\lambda_n - \mu_n) + t^2 \sum_{n=1}^{+\infty} \mu_n  (\lambda_n - \mu_n) + \mathcal O(|b| t^3),
\end{equation}
since $\lambda_n - \mu_n = \mathcal O(\frac{|b|}{n^4})$ and $\mu_n^2 = \mathcal O(n^2)$.

\vspace{0.3cm}\noindent
The expression $(3)$ is easier to study. As previously, we have:
\begin{eqnarray} \label{Asympt3}
	(3) & =& \frac{t^2}{2} \sum_{n=1}^{+\infty} e^{-t\mu_n} (\lambda_n - \mu_n)^2  \nonumber \\
	     &=& \frac{t^2}{2} \sum_{n=1}^{+\infty}  (\lambda_n - \mu_n)^2 + \mathcal O(b^2t^3). 
\end{eqnarray}

\vspace{0.5cm}\noindent
As a conclusion, using (\ref{Asympt1}), (\ref{Asympt2}) and (\ref{Asympt3}), we get immediately:
\begin{eqnarray}\label{tracepositive}
\sum_{n=1}^{+\infty} (e^{-t\lambda_n (b)} -e^{-tn} ) &=& b  + bt \log t + \left( C(b) - \sum_{n=1}^{+\infty} (\lambda_n(b) - \mu_n(b))\ \right) t + b t^2 \log t \nonumber \\ 
                                     & +&  \left( D (b) +  \half \sum_{n=1}^{+\infty}  (\lambda_n^2(b) - \mu_n^2(b) ) \right)  t^2 
                                     + \frac{1}{4} (2b-b^2) \  t^3 \log t + \mathcal O(|b| t^3).
\end{eqnarray}

\vspace{0.3cm}\noindent
It follows that the third term of the (RHS)  of (\ref{trace}) is given by
\begin{eqnarray}\label{tracenegative}
	 \sum_{n=1}^{+\infty} (e^{-t\lambda_n (-b)} -e^{-tn} ) &=& - b  - bt \log t + \left( C(-b) - \sum_{n=1}^{+\infty} (\lambda_n(-b) - \mu_n(-b))\ \right) t  \nonumber \\ 
	& &  - b t^2 \log t + \left( D (-b) +  \half \sum_{n=1}^{+\infty}  (\lambda_n^2(-b) - \mu_n^2 (-b)) \right)  t^2  \nonumber \\
	& & - \frac{1}{4} (2b+b^2) \  t^3 \log t + \mathcal O(|b| t^3).
\end{eqnarray}

\vspace{0.3cm}
\noindent
Thus, thanks to (\ref{trace0}), (\ref{tracepositive}) and (\ref{tracenegative}), we have obtained:
\begin{eqnarray}\label{tracetotale}
{\rm Tr} (e^{-t\Lambda(b)} - e^{-t\Lambda(0)}) &=& \left( C(b)+ C(-b) -\lambda_0(b)- \sum_{n=1}^{+\infty} (\lambda_n(b) - \mu_n(b))\ - \sum_{n=1}^{+\infty} (\lambda_n(-b) - \mu_n(-b)) \right) t \nonumber \\ 
	&+&  \left( D(b)+ D (-b) + \frac{\lambda_0^2}{2} +  \half \sum_{n=1}^{+\infty}  (\lambda_n^2(b) - \mu_n^2(b) ) +  \half \sum_{n=1}^{+\infty}  (\lambda_n^2(-b) - \mu_n^2(-b) ) \right)  t^2 \nonumber \\
	&-& \frac{b^2}{2} t^3 \log t + \mathcal O(|b| t^3).
\end{eqnarray}

\vspace{0.3cm}
\noindent
Using the expressions of $C(b)$ and $D(b)$, we have proved the following result:

\vspace{0.3cm}
\noindent
\begin{thm}\label{tracerelative1}
	Let $L$ be a fixed positive real and $b \in [-L,L]$.  When $t \to 0^+$, we get the asymptotic expansion of the relative Steklov heat trace operator, (and uniformly with respect to small $b$),
	\begin{eqnarray}\label{tracerelative}
	{\rm Tr}(e^{-t\Lambda(b)} - e^{-t\Lambda(0)}) &=&  \left(\big(1-\frac{\pi^2}{3} +7 \zeta(3)\big)\, b^2 -\lambda_0(b)- \sum_{n=1}^{+\infty} (\lambda_n(b) - \mu_n(b))   - \sum_{n=1}^{+\infty} (\lambda_n(-b) - \mu_n(-b)) \right) t \nonumber \\ 
	& & + \left( \frac{1}{3780} \left(    (5670-2520\, \pi^2 +126 \,\pi^4 +4 \, \pi^6 -15120\, \zeta (3)    -7560\, \zeta(5))\, b^2 \right. \right. \nonumber \\
	& &  \left. +(126 \, \pi^4 +57 \,\pi^6 -7560 \,\zeta(3) -34020 \, \zeta(5)) \, b^4 +4 \, \pi^6 \, b^6 \right)  \nonumber \\
	& &  \left. + \frac{\lambda_0^2(b)}{2} +  \half \sum_{n=1}^{+\infty}  (\lambda_n^2(b) - \mu_n^2(b) ) +  \half \sum_{n=1}^{+\infty}  (\lambda_n^2(-b) - \mu_n^2(-b) ) \right)  t^2\nonumber \\
	& & - \frac{b^2}{2} t^3 \log t + \mathcal O(|b|\,t^3).
\end{eqnarray}
\end{thm}

\begin{rem}\label{parity}
Using Mathematica, we are able to compute  the asymptotics of  ${\rm Tr} (e^{-t\Lambda(b)} - e^{-t\Lambda(0)})$ up to the order  $O(t^4)$. The coefficient of $t^3$  is too complicated to be written here, but essentially, this term is in the same form as the coefficient of $t^2$ in Theorem \ref{tracerelative1}, (it is even with respect to $b$), and we can see that the error term in Theorem \ref{tracerelative1} is actually equal to  $\mathcal O(b^2 t^3)$. In contrast, the next logarithmic term  is very simple and is given by  $-b^2 t^4 \log t$. 
\end{rem}

\vspace{0.2cm}
\par\noindent
As a by-product, recalling that ${\rm Tr}(e^{-t\Lambda(0)}) = \coth (\frac{t}{2})$,  we get :
\begin{eqnarray}\label{tracetotale}
	{\rm Tr}(e^{-t\Lambda(b)}) &=& \frac{2}{t}  + \left( \frac{1}{6}+(1-\frac{\pi^2}{3} +7 \zeta(3))b^2 -\lambda_0(b)- \sum_{n=1}^{+\infty} (\lambda_n(b) - \mu_n(b))   - \sum_{n=1}^{+\infty} (\lambda_n(-b) - \mu_n(-b)) \right) t \nonumber \\ 
	& & + \left( \frac{1}{3780} \left(    (5670-2520 \, \pi^2 +126 \, \pi^4 +4 \,\pi^6 -15120 \,\zeta (3)    -7560 \,\zeta(5))\,  b^2 \right. \right. \nonumber \\
	& &  \left. +(126\, \pi^4 +57 \, \pi^6 -7560\,  \zeta(3) -34020 \, \zeta(5))\, b^4 +4 \,\pi^6 \, b^6 \right)  \nonumber \\
	& &  \left. + \frac{\lambda_0^2(b)}{2} +  \half \sum_{n=1}^{+\infty}  (\lambda_n^2(b) - \mu_n^2(b) ) +  \half \sum_{n=1}^{+\infty}  (\lambda_n^2(-b) - \mu_n^2(-b) ) \right)  t^2\nonumber \\
	& & - \frac{b^2}{2} t^3 \log t + \mathcal O(t^3).
\end{eqnarray}

\par \noindent
According to us, it is the first time that an explicit logarithmic term is computed   for the magnetic Steklov heat trace.  For a general pseudo-differential operator of degree $1$, it is known that the presence of logarithmic terms is "generic", (see  Gilkey-Grubb \cite{GiGr1998}). We do not know if such a generic result holds inside the class of D-to-N operators.

\subsection{Some remarks in the regime $ b \to 0$.}

The coefficients of order $\mathcal O(t)$ and $\mathcal O(t^2)$ appearing in Theorem \ref{tracerelative1} are rather cumbersome to study. 
So, in this subsection, we give the asymptotics of these coefficients when the constant magnetic field $b \to 0$. 

\vspace{0.3cm}\noindent
We set :
\begin{equation}\label{coefE}
	E(b) = (1-\frac{\pi^2}{3} +7 \zeta(3))b^2 -\lambda_0 (b)- \sum_{n=1}^{+\infty} (\lambda_n (b) - \mu_n (b))   - \sum_{n=1}^{+\infty} (\lambda_n(-b)- \mu_n (-b)) .
\end{equation}
Using  (\ref{asymptMimp}), we get for $n \geq 1$,
\begin{equation}
\lambda_n (b) - \mu_n (b) = - \frac{1}{n^3 (n+1)} b + \frac{18n^2+31n+14}{2n^3 (n+1)^2(n+2)} b^2 +\mathcal O\big(\frac{|b|^3}{n^4}\big ).
\end{equation}
Thus,  we obtain:
\begin{equation}
	\sum_{n=1}^{+\infty} (\lambda_n (b) - \mu_n (b)) = (1 - \frac{\pi^2}{6} +\zeta(3))\,b + \frac{1}{4}(5-\pi^2 +14 \zeta (3))\, b^2 + \mathcal O(|b|^3).
\end{equation}	

\vspace{0,5cm}\noindent
Then, we deduce easily from (\ref{coefE}) that:
\begin{equation}
	E(b) =  ( \frac{\pi^2}{6} - \frac{7}{4}) \ b^2 + \mathcal O(|b|^3).
\end{equation}

\begin{rem}
Numerically, $\frac{\pi^2}{6} - \frac{7}{4} \simeq -0.105$. Thus, for small enough $b$, we see that $E(b) \leq 0$, and in particular the diamagnetic property holds (in a weak sense).

\vspace{0.3cm}\noindent
In the same way, if we set 
\begin{eqnarray}\label{coefF}
	F(b)= 	& & \frac{1}{3780} \left(    (5670-2520 \, \pi^2 +126\, \pi^4 +4 \, \pi^6 -15120 \,\zeta (3)    -7560 \, \zeta(5))\, b^2  \right. \nonumber \\
	& &  \left. +(126 \, \pi^4 +57 \,\pi^6 -7560 \,\zeta(3) -34020 \, \zeta(5)) \, b^4 +4 \pi^6 \, b^6 \right)  \nonumber \\
	& &  + \frac{\lambda_0^2}{2}(b) +  \half \sum_{n=1}^{+\infty}  (\lambda_n^2(b) - \mu_n^2(b) ) +  \half \sum_{n=1}^{+\infty}  (\lambda_n^2(-b) - \mu_n^2(-b) ),
\end{eqnarray}
we find  with Mathematica that 
\begin{equation}
	F(b) = \frac{2\pi^2 -21}{6}\  b^2 + \mathcal O(|b|^3).
\end{equation}
As previously, since $ \frac{2\pi^2 -21}{6}\simeq -0.210$, we see that the diamagnetic property also holds for this coefficient.
\end{rem}

\vspace{0.5cm}

\noindent \footnotesize{
	
	\noindent Laboratoire de Math\'ematiques Jean Leray, UMR CNRS 6629. Nantes Universit\'e, 2 Rue de la Houssini\`ere, BP 92208, F-44322 Nantes Cedex 03 \\
	\emph{Email adress}: Bernard.Helffer@univ-nantes.fr \\

	\noindent Laboratoire de Math\'ematiques Jean Leray,  UMR CNRS 6629. Nantes Universit\'e, 2 Rue de la Houssini\`ere, BP 92208, F-44322 Nantes Cedex 03 \\
	\emph{Email adress}: francois.nicoleau@univ-nantes.fr \\

\bibliographystyle{acm}
\bibliography{BiblioMagnetique}

\end{document}